\documentclass{tac}

\usepackage{amsmath}
\usepackage{amssymb}
\usepackage{xy}

\input diagxy

\usepackage[colorlinks=true]{hyperref}
\hypersetup{allcolors=[rgb]{0.1,0.1,0.4}}

\author{David Neal Broodryk}

\address{Department of Mathematics and Applied Mathematics, University of Cape Town\\
Rondebosch 7701\\[5pt]
}

\title{Characterization of coextensive varieties of universal algebras II}

\copyrightyear{2020}

\keywords{Coextensivity, Universal Algebra, Syntactic Characterization}
\amsclass{18A30, 08B05}

\eaddress{BRDDAV020@myuct.ac.za}

\newtheorem{thm}{Theorem}

\newtheoremrm{rem}{Remark}
\newtheorem{corollary}{Corollary}

\let\pf\proof
\let\epf\endproof

\newcommand{\N}{\mathbb{N}}

\newcommand{\C}{\mathcal{C}}
\newcommand{\CSR}{\textbf{CSemiRing}}
\newcommand{\Fo}{F(\emptyset)}
\newbox\langlebox
\setbox\langlebox=\hbox{\xy \POS(75,75)\ar@{-} (0,75) \ar@{-} (75,0)\endxy}
\def\langle{\copy\langlebox}
\newbox\ranglebox
\setbox\ranglebox=\hbox{\xy \POS(0,75)\ar@{-} (0,0) \ar@{-} (75,75)\endxy}
\def\rangle{\copy\ranglebox}

\begin{document}

\maketitle

\begin{abstract}
A coextensive category can be defined as a category $\C$ with finite products such that for each pair $X,Y$ of objects in $\C$, the canonical functor $\times\colon X/\C \times Y/\C \to (X \times Y)/\C$ is an equivalence. We give a syntactical characterization of coextensive varieties of universal algebras. This paper is an updated version of the pre-print \cite{COEX}.
\end{abstract}

\section{Coextensivity and equivalent conditions}

We say that a Category $\C$ is Co-extensive if for any objects $X,Y$ in $\C$, the canonical functor $\times \colon X/\C \times Y/\C \to (X \times Y)/\C$ is an equivalence. Here $\times \colon X/\C \times Y/\C \to (X \times Y)/\C$ is the functor that maps a pair of morphisms $(f,g)$ to their product $f\times g$. The left adjoint $L\colon  (X \times Y)/\C \to X/\C \times Y/\C$ sends a morphism $h: X\times Y \to Z$ to its pushouts across $\pi_X$ and $\pi_Y$.
\begin{equation}
\bfig
\node 00(0,0)[X]
\node 01(600,0)[X\times Y]
\node 02(1200,0)[Y]
\node 03(1800,0)[X]
\node 04(2400,0)[X\times Y]
\node 05(3000,0)[Y]
\node 10(0,-600)[A]
\node 11(600,-600)[A\times B]
\node 12(1200,-600)[B]
\node 13(1800,-600)[A']
\node 14(2400,-600)[Z]
\node 15(3000,-600)[B']
\node (1875,-525)[\langle]
\node (2925,-525)[\rangle]
\arrow[01`00;\pi_X]
\arrow[01`02;\pi_Y]
\arrow[04`03;\pi_X]
\arrow[04`05;\pi_Y]
\arrow[00`10;f]
\arrow[01`11;f\times g]
\arrow[02`12;g]
\arrow[03`13;L(h)_1]
\arrow[04`14;h]
\arrow[05`15;L(h)_2]
\arrow[11`10;\pi_A]
\arrow[11`12;\pi_B]
\arrow[14`13;p_1']
\arrow[14`15;p_2']
\efig
\end{equation}
The existence of the left adjoint $L$ gives rise to multiple conditions equivalent to coextensivity. Let $\eta$ and $\varepsilon$ be the unit and counit of this adjunction. Then $\C$ is co-extensive if and only if $\varepsilon$ and $\eta$ are natural isomorphisms as in the diagrams:
\begin{equation}
\bfig
\node 00(0,0)[X]
\node 01(600,0)[X\times Y]
\node 02(1200,0)[Y]
\node 03(1800,0)[X]
\node 04(2400,0)[X\times Y]
\node 05(3000,0)[Y]
\node 10(0,-600)[A]
\node 11(600,-600)[A\times B]
\node 12(1200,-600)[B]
\node 13(1800,-600)[A']
\node 14(2400,-600)[Z]
\node 15(3000,-600)[B']
\node 20(-400,-1000)[C]
\node 21(1600,-1000)[D]
\node 22(2400,-1000)[A'\times B']
\node 23(3400,0)[]
\node (1875,-525)[\langle]
\node (2925,-525)[\rangle]
\arrow[01`00;\pi_X]
\arrow[01`02;\pi_Y]
\arrow[04`03;\pi_X]
\arrow[04`05;\pi_Y]
\arrow[00`10;f]
\arrow[01`11;f\times g]
\arrow[02`12;g]
\arrow[03`13;L(h)_1]
\arrow[04`14;h]
\arrow[05`15;L(h)_2]
\arrow[11`10;\pi_A]
\arrow[11`12;\pi_B]
\arrow[14`13;p_1']
\arrow[14`15;p_2']
\arrow|m|[00`20;L(f\times g)_1]
\arrow|m|[02`21;L(f\times g)_2]
\arrow|m|[20`10;\varepsilon_f]
\arrow|m|[21`12;\varepsilon_g]
\arrow|b|[11`20;p_1]
\arrow|b|[11`21;p_2]
\arrow|b|[22`13;\pi_A']
\arrow|b|[22`15;\pi_B']
\arrow|m|[14`22;\eta_h]
\efig
\end{equation}

Equivalently, since both products and pushouts are defined up to isomorphism, $\C$ is co-extensive when, for any objects $X$ and $Y$, and any commutative diagram:

\begin{equation}
\bfig
\node 03(1800,0)[X]
\node 04(2400,0)[X\times Y]
\node 05(3000,0)[Y]
\node 13(1800,-600)[A']
\node 14(2400,-600)[Z]
\node 15(3000,-600)[B']
\arrow[04`03;\pi_X]
\arrow[04`05;\pi_Y]
\arrow[03`13;f]
\arrow[04`14;h]
\arrow[05`15;g]
\arrow[14`13;p_1']
\arrow[14`15;p_2']
\efig
\end{equation}

the bottom row of the diagram is a product diagram if and only if, both squares are pushouts \cite{CLW}. In this case we say that $\C$ satisfies the coextensivity condition for any $X$ and $Y$. When $\C$ has an initial object $0$ we can simplify this condition even further.

\proposition\label{proposition_initial} If a category with products has an initial object $0$, then it is coextensive if and only if it satisfies the coextensivity condition for $0$ and $0$.
\endthm

\pf This is simply the dual condition of proposition 4.1 of \cite{CLW}. When $\C$ is coextensive then it satisfies the coextensivity condition for any $X$ and $Y$, including $0$ and $0$. On the other hand, consider the diagram

\begin{equation}
\bfig
 \node 00(0,0)[0]
 \node 01(600,0)[0\times0]
 \node 02(1200,0)[0]
 \node 10(0,-600)[X]
 \node 11(600,-600)[X\times Y]
 \node 12(1200,-600)[Y]
 \node 20(0,-1200)[A']
 \node 21(600,-1200)[S]
 \node 22(1200,-1200)[B']
 \node(300,-300)[(1)]
 \node(900,-300)[(3)]
 \node(300,-900)[(2)]
 \node(900,-900)[(4)]
 \arrow[01`00;\pi_1]
 \arrow[01`02;\pi_2]
 \arrow[11`10;\pi_X]
 \arrow[11`12;\pi_Y]
 \arrow[21`20;p_1']
 \arrow[21`22;p_2']
 \arrow[00`10;!_X]
 \arrow|m|[01`11;!_X \times !_Y]
 \arrow|r|[02`12;!_Y]
 \arrow[10`20;f]
 \arrow[11`21;h]
 \arrow|r|[12`22;g]
\efig
\end{equation}

where $!_X$ and $!_Y$ are the unique morphisms from $0$ to $X$ and $Y$ respectively. If $\C$ satisfies the coextensivity condition for $0$ and $0$ then (1) and (3) are pushouts for all $X,Y$ in $\C$. Now, by the pasting law, (2) and (4) are pushouts $\iff$ (1)+(2) and (3)+(4) are pushouts $\iff$ the bottom row is a product diagram. Therefore $\C$ is coextensive.

\epf

There are also multiple ways of splitting coextensivity into smaller conditions. We say that $\C$ is:

\begin{description}
\item[Left coextensive] when the co-unit $\varepsilon$ is a natural isomorphism.
\item[Right coextensive] when the unit $\eta$ is a natural isomorphism.
\end{description}

Equivalently $\C$ is left coextensive when, in the previous diagram, both squares are pushouts whenever the bottom row is a product diagram. Similarly $\C$ is right coextensive when, in the previous diagram, the bottom row is a product diagram whenever both squares are pushouts. Clearly $\C$ is coextensive if and only if $\C$ is left and right coextensive. Alternatively we say that $\C$ has:

\begin{description}
\item[co-universal products] when $\C$ is right-coextensive ($\eta$ is a natural isomorphism)
\item[co-disjoint products] when all product projections are epimorphisms and the following diagram is a pushout for all $X$ and $Y$ in $\C$:

\begin{equation}
\bfig
\node 00(0,0)[X\times Y]
\node 01(600,0)[Y]
\node 10(0,-600)[X]
\node 11(600,-600)[1]
\node (525,-525)[\rangle]
\arrow[00`01;\pi_2]
\arrow[00`10;\pi_1]
\arrow[10`11;]
\arrow[01`11;]
\efig
\end{equation}
\end{description}

Even though $\C$ having co-disjoint products is a weaker condition than $\C$ being left-coextensive, it is still the case that $\C$ is coextensive if and only if $\C$ has co-disjoint and co-universal products \cite{CLW}.

\section{Coextensivity in varieties of universal algebras}

The motivating example of a coextensive variety of universal algebras is the variety $\CSR$ of commutative semi-rings. Let $h: X\times Y \to Z$ be a homomorphism in in this variety. Then $Z$ contains some elements $e_1 = h(1,0)$ and $e_2 = h(0,1)$. Since $h$ is a homomorphism these elements satisfy the identities:

\begin{itemize}
	\item  $e_1+e_2 = 1$
	\item $e_1\cdot e_2 = 0$
	\item $e_1\cdot e_1 = e_1$
	\item $e_2\cdot e_2 = e_2$
\end{itemize}

We can now define the two cosets $e_1Z = \{e_1z | z \in Z\}$, $e_2Z = \{e_2z | z \in Z\}$. These sets can be made into commutative semirings using the same addition, multiplication and $0$ as $Z$. The unit is the only thing that changes, specifically for $i \in \{1,2\}$ and any $e_iz \in e_iZ$ we have $e_i \cdot e_iz = (e_ie_i)z=e_iz$ since both $e_1$ and $e_2$ are idempotent. So $e_1$ and $e_2$ are the units of $e_1Z$ and $e_2Z$ respectively. Note that these operations are defined in just the right way so that the maps $p_i:Z \to e_iZ$ defined by $p_i(z) = e_iz$ are homomorphisms, and in fact both squares in the following diagram are pushouts, where $h_1(x) = h((x,0))$ and $h_2(y) = h((0,y))$.

\begin{equation}
\bfig
\node 03(1800,0)[X]
\node 04(2400,0)[X\times Y]
\node 05(3000,0)[Y]
\node 13(1800,-600)[e_1Z]
\node 14(2400,-600)[Z]
\node 15(3000,-600)[e_2Z]
\arrow[04`03;\pi_X]
\arrow[04`05;\pi_Y]
\arrow[03`13;h_1]
\arrow[04`14;h]
\arrow[05`15;h_2]
\arrow[14`13;p_1]
\arrow[14`15;p_2]
\efig
\end{equation}

It remains to show that the bottom row is a product diagram. To show that $p_1 \times p_2$ is an isomorphism, consider the map $f: e_1Z \times e_2Z \to Z$ defined by $f(e_1z_1,e_2z_2) = e_1z_1+e_2z_2$ for all $z_1,z_2 \in Z$. Then, by the above identities, we have:

\begin{align*}
	f\circ (p_1\times p_2)(z) &= e_1z + e_2z = (e_1+e_2)z = 1z = z \\
	(p_1\times p_2) \circ f(e_1z_1,e_2z_2) &= (e_1e_1z_1+e_1e_2z_2,e_2e_1z_1+e_2e_2z_2) = (e_1z_1,e_2z_2)
\end{align*}

and so $\CSR$ has co-universal products. That $\CSR$ has co-disjoint products follows as a consequence of the following proposition.  

\proposition\label{proposition_disjoint} A variety $\C$ has co-disjoint products if and only if $\C$ has at least one constant term.
\endthm

\pf Note that when $\C$ is a variety, the terminal object $1$ is any singleton considered as an algebra, and diagram $5$ is a pushout as long as neither $X$ nor $Y$ is empty. On the other hand if both $X$ and $Y$ are empty the diagram is never a pushout. Therefore $\C$ has co-disjoint products  exactly when $\emptyset \notin \C$, which is to say that $\C$ has at least one constant term.
\epf

Left coextensivity, on the other hand, has a more involved characterisation which we describe in \cite{DB}. It remains to characterise right coextensivity. Consider the forgetful functor $U: \C \to Sets$ taking each algebra to its underlying set, and its left adjoint, the free functor $F: Sets \to \C$ taking each set $X$ to the free algebra $F(X)$ in $\C$ generated on $X$. Since $\emptyset$ is the initial object in $Sets$ we have that $\Fo$ is the initial object in $\C$. This allows us to consider only the coextensivity condition for $\Fo$ and $\Fo$ as follows from proposition 1.1. It now becomes useful to consider the following diagram in which both squares are pushouts:

\begin{equation}
\bfig
 \node 00(0,0)[\Fo]
 \node 01(800,0)[\Fo \times \Fo]
 \node 02(1600,0)[\Fo]
 \node 10(0,-600)[F(X)]
 \node 11(800,-600)[F(X) + \Fo^2]
 \node 12(1600,-600)[F(X)]
 \node 21(800,-1200)[F(X)^2]
 \node (100,-500)[\langle]
 \node (1500,-500)[\rangle]
 
 \arrow[01`00; \pi_1]
 \arrow[01`02; \pi_2]
 \arrow[00`10;]
 \arrow[01`11;i]
 \arrow[02`12;]
 \arrow[11`10;{[}id,\pi_1{]}]
 \arrow[11`12;{[}id,\pi_2{]}]
 \arrow[11`21;\eta_X]
 \arrow[21`10;p_1]
 \arrow[21`12;p_2]
\efig
\end{equation}

Here $F(X) + \Fo^2$ is the coproduct of $F(X)$ with $\Fo \times \Fo$ and $i$ is the corresponding coproduct inclusion. For ease of notation we write the component of $\eta$ at $i$ as $\eta_X$ instead of as $\eta_i$. Informally, $F(X) + \Fo^2$ can be thought of as the `least restricted' algebra equipped with a morphism from some product, and using this diagram we can now prove the following result: 

\proposition\label{proposition}  Any coextensive variety contains, for some $k\geq 1$, a $(k+2)$-ary term $t$ and constants $e_1, \dots,e_k,e_1', \dots,e_k' \in \Fo$ such that the identities hold:
 \begin{align*}
    t(x, y, e_1, \dots, e_k) &= x \\
    t(x, y, e_1', \dots, e_k') &= y
\end{align*}
\endthm

\pf $F(X) + \Fo^2$ is generated by $X$ and $\Fo^2$, so the elements of $F(X) + \Fo^2$ are all of the form $t(x_1,\dots,x_n,(e_1,e_1'),\dots,(e_k,e_k'))$ for some $n,k \in \N$ where $t$ is some $(k+n)$-ary operation, $x_1,\dots,x_n \in X$, and $e_1,\dots,e_k,e_1',\dots,e_k' \in \Fo$. $\eta_X$ is the unique map $([id,\pi_1],[id,\pi_2]) = [\Delta_{F(X)}, \subseteq_{\Fo^2}]$ such that the diagram commutes. Therefore: $$\eta_X(t(x_1,\dots,x_n,(e_1,e_1'),\dots,(e_k,e_k'))) = (t(x_1,\dots,x_n,e_1,\dots,e_k), t(x_1,\dots,x_n,e_1',\dots,e_k'))$$
Now, let $x,y \in X$, then $(x,y) \in F(X)^2$, so since $\eta_X$ is an isomorphism there must exist some $t(x,y,(e_1,e_1'),\dots,(e_k,e_k'))$ such that $$\eta_X(t(x,y,(e_1,e_1'),\dots,(e_k,e_k'))) = (t(x,y,e_1,\dots,e_k),t(x,y,e_1',\dots,e_k')) = (x,y)$$

which gives the desired identities. Since $\C$ is co-extensive, $\C$ has co-disjoint products and so $\C$ has at least one constant term. We can therefore assume that $k\geq1$

\epf

We say that $t$ is a \textbf{diagonalising term}, and when such a term exists we can  define the map $\delta: F(X)^2 \to F(X) + \Fo^2$ as $\delta(x,y) = t(x,y,(e_1,e_1'),\dots,(e_n,e_n'))$. Note that $\eta_X\delta((x,y)) = (x,y)$, so $\delta$ is a right inverse of $\eta_X$ considered as a map, and in particular $\eta_X$ is an isomorphism if and only if $\delta$ is. This turns out to be especially important because of the next two results.

\proposition\label{proposition_diag_conditions} Let $\C$ be a variety with a diagonalising term $t$. then $\C$ satisfies the following conditions:
\begin{enumerate}
	\item Any reflexive homomorphic relation $R$ on some product $A \times B \in \C$ is of the form $R_A \times R_B$ for some reflexive homomorphic relations $R_A$ on $A$ and $R_B$ on $B$.
	\item Any congruence $E$ on some product $A \times B \in \C$ is of the form $E_A \times E_B$ for some congruences $E_A$ on $A$ and $E_B$ on $B$.
	\item $\C$ satisfies the coextensivity condition for any surjective homomorphism $q$.
\end{enumerate}

\endthm

\pf
(1) Let $R$ be a reflexive homomorphic relation on $A\times B$, then $R \subseteq (A\times B)^2 \simeq A^2 \times B^2$. Let $\pi_1 : A^2 \times B^2 \to A^2$ and $\pi_2: A^2 \times B^2 \to B^2$ be product projections. Then the images of $R$ under $\pi_1$ and $\pi_2$, $R_A = \pi_1(R)$, and $R_B = \pi_2(R)$, are reflexive homomorphic relations on $A$ and $B$ respectively. Concretely, these relations are of the form:
\begin{align*}
	R_A &= \{(a,c) \in A^2 | ((a,b),(c,d)) \in R \text{ for some } b,d \in B\} \\
	R_B &= \{(b,d) \in A^2 | ((a,b),(c,d)) \in R \text{ for some } a,c \in A\}
\end{align*}
Therefore, their product $R_A \times R_B$ is of the form:
\begin{align*}
	((a,b),(c,d)) \in R_A \times R_B &\iff (a,c) \in R_A \text{ and } (b,d) \in R_B \\
	&\iff ((a,b'),(c,d')) \in R \text{ and } ((a',b),(c',d)) \in R
\end{align*}
For some $a',c' \in A$ and $b',d' \in B$. Clearly $R \subseteq R_A \times R_B$. On the other hand, whenever $(a,b)\simeq_{R_A \times R_B}(c,d)$ we have that $((a,b'),(c,d')) \in R$ and $((a',b),(c',d)) \in R$. Then, since $R$ is a reflexive homomorphic relation, we have that:
$$((a,b),(c,d)) = t(((a,b'),(c,d')),((a',b),(c',d)),((e_1,e_1'),(e_1,e_1')),\dots,((e_k,e_k'),(e_k,e_k')) \in R$$
and so $R = R_A \times R_B$ as desired.

\vspace{1em}

(2) Note that when $E$ is a congruence, then $E_A = \pi_1(E)$ and $E_B = \pi_2(E)$ are immediately reflexive, symmetric and homomorphic. In general, they may not be transitive, but if $\C$ has a diagonalising term then $E = E_A \times E_B$ by (1), and so $a_1\simeq_{E_A}a_2 \simeq_{E_A}a_3$ implies $(a_1,b) \simeq_E (a_2,b) \simeq_E (a_3,b)$ for any $b \in B$, and so $(a_1,b) \simeq_E (a_3,b)$ by transitivity of $E$, which implies $a_1 \simeq_{E_A} a_3$ as desired. 

\vspace{1em}

(3) Consider some surjective homomorphism $q: A\times B \to C$ for some algebras $A,B$ and $C$ in $\C$. Note that $C \simeq (A\times B)/E$ where $E$ is the congruence on $A\times B$ generated by $q$. We can now construct the following diagram in which both squares are pushouts:

\begin{equation}
\bfig
 \node 00(0,0)[A]
 \node 01(600,0)[A \times B]
 \node 02(1200,0)[B]
 \node 10(0,-600)[A/E_A]
 \node 11(600,-600)[(A \times B) / E]
 \node 12(1200,-600)[B/E_B]
 \node (100,-500)[\langle]
 \node (1100,-500)[\rangle]
 
 \arrow[01`00; \pi_1]
 \arrow[01`02; \pi_2]
 \arrow[00`10;]
 \arrow[01`11;q]
 \arrow[02`12;]
 \arrow[11`10;]
 \arrow[11`12;]
\efig
\end{equation}

Note that $(A/E_A) \times (B/E_B) \simeq (A\times B)/(E_A \times E_B) \simeq (A\times B)/E$ by the previous result, and so the bottom row is a product diagram. Therefore $\C$ satisfies the coextensivity condition for $q$ as desired.
\epf

\proposition\label{proposition_sufficient} A variety $\C$ with a diagonalising term $t$ is coextensive if and only if $\eta_X$ is an isomorphism for all sets $X$.
\endthm

\pf
If $\C$ is coextensive, then $\eta_X$ is immediately an isomorphism for all sets $X$. On the other hand, given a morphism $i: \Fo^2 \to A$ for some algebra $A$, let $U(A)$ be the underlying set of $A$ and $\varepsilon_A: FU(A) \to A$ the canonical morphism from the free algebra on $U(A)$ to $A$. Then $i$ can be factorised as the composition of the morphisms:

\begin{equation}
\bfig
\node 0(0,0)[\Fo^2]
\node 1(1000,0)[FU(A)+\Fo^2]
\node 2(2000,0)[A]
\arrow[0`1;\subset]
\arrow[1`2;{[}\varepsilon_A,i{]}]
\efig
\end{equation}

Furthermore, since $\varepsilon_A$ is surjective, so is ${[}\varepsilon_A,i{]}$. We can now construct the following diagram in which all four squares are pushouts.

\begin{equation}
\bfig
 \node 00(0,0)[\Fo]
 \node 01(800,0)[\Fo \times \Fo]
 \node 02(1600,0)[\Fo]
 \node 10(0,-600)[FU(A)]
 \node 11(800,-600)[FU(A) + \Fo^2]
 \node 12(1600,-600)[FU(A)]
 \node 20(0,-1200)[A_1]
 \node 21(800,-1200)[A]
 \node 22(1600,-1200)[A_2]
 \node (100,-500)[\langle]
 \node (100,-1100)[\langle]
 \node (1500,-500)[\rangle]
 \node (1500,-1100)[\rangle]
 \arrow[01`00; \pi_1]
 \arrow[01`02; \pi_2]
 \arrow[00`10;]
 \arrow[01`11;]
 \arrow[02`12;]
 \arrow[11`10; {[}id,\pi_1{]}]
 \arrow[11`12; {[}id,\pi_2{]}]
 \arrow[10`20;]
 \arrow[11`21;{[}\varepsilon_A,i{]}]
 \arrow[12`22;]
 \arrow[21`20;]
 \arrow[21`22;]
\efig
\end{equation}

Note that if $\eta_X$ is an isomorphism for all sets $X$, then $\eta_{U(A)}$ is an isomorphism and the middle row is a product diagram. Then since ${[}\varepsilon_A,i{]}$ is surjective, the bottom row is also a product diagram as desired by the previous proposition. 

\epf

This result establishes the sufficiency of diagram $6$ as the only diagram we need consider, since $\C$ will be coextensive if and only if $\eta_X$ is an isomorphism for any set $X$. Recall that whenever $\C$ has a diagonalising term we can define the map $\delta: F(X)^2 \to F(X) + \Fo^2$ as $\delta(x,y) = t(x,y,(e_1,e_1'),\dots,(e_n,e_n'))$, and so $\eta_X \delta ((x,y)) = (x,y)$. Our next step is to convert our previous result into one concerning $\delta$ instead of $\eta_X$.

\proposition\label{proposition_delta} 
A variety $\C$ is coextensive if and only if for any set $X$:
\begin{enumerate}
    \item $\C$ has a diagonalising term $t$
    \item $\delta(x,x) = x$ for all $x \in X$     
    \item $\delta:F(X)^2 \to F(X) + \Fo^2$ is a homomorphism
\end{enumerate}
\endthm

\pf It follows from propositions $2.2$ and $2.4$ that a variety $\C$ is coextensive if and only if $\C$ has a diagonalising term $t$ and $\eta_X$ is an isomorphism for any set $X$. Since $\delta$ is a right inverse of $\eta_X$ considered as a map, $\eta_X$ is an isomorphism if and only if $\delta$ is an isomorphism. In this case $\delta$ is immediately a homomorphism, and $\delta(x,x) = \delta\eta_X(x) = x$ as desired.

\vspace{1em}

On the other hand, assume that $\delta$ is a homomorphism and that $\delta(x,x)=x$. Since $\delta$ has a left inverse, $\delta$ is already injective. Thus, it remains to show that $\delta$ is surjective. For any $(w,w') \in \Fo^2$ we have:
\begin{align*}
    \delta(w,w') &= t((w,w),(w',w'),(e_1,e_1'),\dots,(e_k,e_k')) \\
    &= (t(w,w',e_1,\dots,e_k),t(w,w',e_1',\dots,e_k')) \\
    &= (w,w')
\end{align*}
Every element of $F(X)+\Fo^2$ is, for some operation $u$, elements $x_1,\dots,x_n \in F(X)$ and $(e_1,e_1'),\dots,(e_m,e_m') \in \Fo^2$, of the form:
\begin{align*}
    u(x_1,\dots,x_n,(e_1,e_1'),\dots,(e_m,e_m')) &= u(\delta(x_1,x_1),\dots,\delta(x_n,x_n),\delta(e_1,e_1'),\dots,\delta(e_m,e_m')) \\
    &= \delta(u(x_1,\dots,x_n,(e_1,e_1'),\dots,(e_m,e_m')))
\end{align*}
Since $\delta$ is a homomorphism satisfying $\delta(x,x)=x$. Therefore $\delta$ is surjective and thus an isomorphism as desired.  
\epf

When converting a categorical condition into a syntactical characterisation it is often useful to first convert the condition into a statement about free algebras, since they obey only the identities that are satisfied by the entire variety. In our case, we have converted coextensivity into a statement about the algebra $F(X) + \Fo^2$. This is not quite a free algebra, as it needs to accept a morphism from $\Fo^2$, so instead we consider it as a quotient of a free algebra. 

\proposition\label{proposition_C} Let $C = \{c_1,\dots,c_k\}$ be a set with $k$ elements. Then $F(X) + \Fo^2 \simeq F(X \cup C)/E$, where $E$ is the congruence generated the relation:
$$R = \{ v(c_1,\dots,c_k) \simeq t(v(e_1,\dots,e_k),v(e_1',\dots,e_k'),c_1,\dots,c_k) \text{ for all $k$-ary operations v} \} $$

\endthm

\pf Let $p$ be the unique morphism $p: F(X \cup C) \to F(X) + \Fo^2$ such that $p(x) = x$ for all $x \in X$ and $p(c_i) = (e_i,e_i')$ for all $i \leq k$. For any $\omega,\omega' \in \Fo$ we have: $$p(t(\omega,\omega',c_1,\dots,c_k)) = t((\omega,\omega),(\omega',\omega'),(e_1,e_1'),\dots,(e_k,e_k')) = (\omega,\omega')$$ For any $x \in F(X)$ we simply have $p(x)=x$. Therefore, since $F(X) + \Fo^2$ is generated by $F(X)$ and $\Fo^2$, $p$ is surjective and so $F(X) + \Fo^2 \simeq F(X \cup C)/E_p$, where $E_p$ is the congruence generated by $p$. It remains to show that $E_p$ is generated by $R$. Note that $$p(v(c_1,\dots,c_k)) = (v(e_1,\dots,e_k),v(e_1',\dots,e_k')) = p(t(v(e_1,\dots,e_k),v(e_1',\dots,e_k'),c_1,\dots,c_k))$$ so $R \subseteq E_p$. On the other hand, for any $v(c_1,\dots,c_k), u(c_1,\dots,c_k) \in F(C)$ we have that $v(c_1,\dots,c_k) \simeq_{E_p} u(c_1,\dots,c_k)$ implies $p(v(c_1,\dots,c_k))=p(u(c_1,\dots,c_k))$, which gives:
\begin{align*}
    (v(e_1,\dots,e_k),v(e_1',\dots,e_k')) = (u(e_1,\dots,e_k),u(e_1',\dots,e_k'))
\end{align*}
So $v(e_1,\dots,e_k) = u(e_1,\dots,e_k)$ and $v(e_1',\dots,e_k') = u(e_1',\dots,e_k')$. Since these identities hold in $\Fo$, they must hold in any algebra in the variety. Let $E$ be the congruence on $F(X \cup C)$ generated by $R$. Then:
\begin{align*}
    v(c_1,\dots,c_k) &\simeq_{E} t(v(e_1,\dots,e_k),v(e_1',\dots,e_k'),c_1,\dots,c_k) \\
    &= t(u(e_1,\dots,e_k),u(e_1',\dots,e_k'),c_1,\dots,c_k) \\
    &\simeq_{E} u(c_1,\dots,c_k)
\end{align*}

and so $E = E_p$ and $F(X) + \Fo^2 \simeq F(X \cup C)/E$ as desired.
\epf

We can now give our characterisation of coextensive varieties of universal algebras.

\thm\label{theorem} Let $\C$ be a variety with a diagonalising term $t$ and constants $e_1,\dots,e_k$ and $e_1',\dots,e_k'$ for some $k\geq 1$. Given two terms $\alpha$ and $\beta$, we write $\alpha \simeq_t \beta$ if $\C$ contains, for some $n,m,m' \in \N$, a sequence $u_1,\dots u_n$ of $(m+m')$-ary terms, and sequences $v_1,\dots,v_{m'}$ and $w_1,\dots,w_{m'}$ of $k$-ary terms, such that the following identities hold in $\C$:
\begin{align*}
	u_1(x_1,\dots,x_m,w_1(c_1,\dots,c_k),\dots,w_{m'}(c_1,\dots,c_k)) &= \alpha \\
    u_i(x_1,\dots,x_m,v_1(c_1,\dots,c_k),\dots,v_{m'}(c_1,\dots,c_k)) &= u_{i+1}(x_1,\dots,x_m,w_1(c_1,\dots,c_k),\dots,w_{m'}(c_1,\dots,c_k)) \\
    u_n(x_1,\dots,x_m,v_1(c_1,\dots,c_k),\dots,v_{m'}(c_1,\dots,c_k)) &= \beta
\end{align*}
And for each $1 \leq j \leq m'$, one of the equalities holds: 
\begin{align*}
	v(c_1,\dots,c_k) &= t(w(e_1,\dots,e_k),w(e_1',\dots,e_k'),(c_1,\dots,c_k)) \\ 
	w(c_1,\dots,c_k) &= t(v(e_1,\dots,e_k),v(e_1',\dots,e_k'),(c_1,\dots,c_k))
\end{align*}

A variety $\C$ is coextensive if and only if it satisfies the following conditions:
\begin{enumerate}
	\item $\C$ contains a diagonalising term $t$
	\item $t(x,x,c_1,\dots,c_k) \simeq_t x$ for any $x,c_1,\dots,c_k$
	\item $s(t(x_1,y_1,c_1,\dots,c_k),\dots,t(x_p,y_p,c_1,\dots,c_k)) \simeq_t t(s(x_1,\dots,x_p),s(y_1,\dots,y_p),c_1,\dots,c_k)$

	for any $p \in \N$, any $p$-ary operation $s$ in $\C$ and any $x_1,\dots,x_p,y_1,\dots,y_p,c_1,\dots,c_k$
\end{enumerate}

\endthm
 
\pf It follows from proposition 2.6 that $\delta(x,x) = x$ if and only if $t(x,x,c_1,\dots,c_k) \simeq_{E} x$, and similarly $\delta$ is a homomorphism if and only if:
\begin{align*}
    \delta(s(x_1,\dots,x_p),s(y_1,\dots,y_p)) &= s(\delta(x_1,y_1),\dots,\delta(x_p,y_p)) \iff \\
    t(s(x_1,\dots,x_p),s(y_1,\dots,y_p),c_1,\dots,c_k) &\simeq_{E} s(t(x_1,y_1,c_1,\dots,c_k),\dots,t(x_p,y_p,c_1,\dots,c_k))
\end{align*}

for all $p\in\N$, $p$-ary operations $s$ and all $x_1,y_1,\dots,x_p,y_p \in X$. It remains to show that given $\alpha,\beta \in F(X\cup C)$, we have $\alpha \simeq_E \beta$ if and only if $\alpha \simeq_t \beta$. $E$ is the congruence generated by $R$, and thus it is the transitive, homomorphic, symmetric, reflexive closure of $R$. Therefore $(\alpha,\beta)\in E$ if and only if $(\alpha,\beta) \in Q^n$ for some $n \in \N$, where $Q$ is the homomorphic, symmetric, reflexive closure of $R$. But $(\alpha,\beta) \in Q^n$ if and only if there exists a sequence $a_0,\dots,a_n \in F(X\cup C)$ such that $a_0 = \alpha$, $a_n = \beta$, and $(a_i,a_{i+1}) \in Q$ for all $i < n$. Finally note that some pair $(a_i,a_{i+1}) \in Q$ if and only if for some operation $u_i$ we have:
\begin{align*}
    a_i &= u_i(x_1,\dots,x_m,v_1(c_1,\dots,c_k),\dots,v_{m'}(c_1,\dots,c_k)) \\
    a_{i+1} &= u_i(x_1,\dots,x_m,w_1(c_1,\dots,c_k),\dots,w_{m'}(c_1,\dots,c_k))
\end{align*}
where either $(v_j(c_1,\dots,c_k),w_j(c_1,\dots,c_k)) \in R$ or $(w_j(c_1,\dots,c_k),v_j(c_1,\dots,c_k)) \in R$ for each $j \leq m'$. This is exactly to say that $\alpha \simeq_t \beta$.
\epf

This concludes the characterisation of coextensive varieties of universal algebras, which is the main result of this paper. Since we also already know the characterisations of left-coextensive varieties, as well as varieties with co-disjoint products, it is natural to wonder about right coextensive varieties that aren't coextensive. In fact there is only one such variety.

\proposition\label{proposition_right} A category is right coextensive without being coextensive if and only if it has no constant terms and satisfies the identity $x=y$. Equivalently, the variety contains only the empty algebra and singleton algebras.
\endthm

\pf Let $\C$ be right coextensive without being coextensive, this is to say that $\C$ has co-universal products but not co-disjoint products. Therefore $\C$ has no constant terms and so the initial object $0$ is the empty algebra. Then for any $X \in \C$ we construct the diagram:

\begin{equation}
\bfig
\node 03(1800,0)[0]
\node 04(2400,0)[0]
\node 05(3000,0)[0]
\node 13(1800,-600)[X]
\node 14(2400,-600)[X]
\node 15(3000,-600)[X]
\arrow[04`03;id_0]
\arrow[04`05;id_0]
\arrow[03`13;!_X]
\arrow[04`14;!_X]
\arrow[05`15;!_X]
\arrow[14`13;id_X]
\arrow[14`15;id_X]
\efig
\end{equation}

In which the top row is a product diagram and both squares are pushouts, so the bottom row is also a product diagram. Then, for any $x,y \in X$, there must exist some $z\in X$ such that $x = id_X(z) = y$, and so $x=y$ for all $x,y \in X$ and $X$ is a singelton as desired. Note that up to isomorphism there are only two objects $0,1 \in \C$ and so it is easy to check in each case that right coextensivity holds.

\epf

\section{Diagonalising terms}

As we have seen above, every coextensive variety $\C$ has a diagonalising term $t$ equipped with some constants $e_1,e_1',\dots,e_k,e_k' \in \Fo$ such that the two identities hold:
 \begin{align*}
    t(x, y, e_1, \dots, e_k) &= x \\
    t(x, y, e_1', \dots, e_k') &= y
\end{align*}
Recall that the motivating example of a coextensive variety is the variety of commutative semmirings. In a similar fashion, the best example of a variety with a diagonalising term is the variety of lattices.

\definition\label{def_lattice} A lattice $L$ can be defined as a set $L$ equipped with two constants $0$ and $1$, and two binary, commutative, and associative operations $\vee,\wedge$, satisfying the four identities:
\begin{itemize}
	\item $a \vee (a \wedge b) = a$
	\item $a \wedge (a \vee b) = a$
	\item $a \vee 0 = a$
	\item $a \wedge 1 = a$
\end{itemize}
\endthm

From these identities it follows that $\vee$ and $\wedge$ are idempotent and that they satisfy the identities:
\begin{align*}
a \vee 1 = 1 \\
a \wedge 0 = 0
\end{align*}
Then it is simple to show the following:

\proposition\label{lattice_diag} The variety \textbf{Lat} of lattices has a diagonalising term, \newline $t(a,b,c,d) = (a\wedge c) \vee (b\wedge d)$, with constants $e_1=1,e_2=0,e_1'=1,e_2'=0$ 
\endthm

\pf For any lattice $L \in Lat$ and any $x,y \in L$ we have \begin{align*}
	t(x,y,e_1,e_2) &= (x\wedge 1) \vee (y\wedge 0) = x \vee 0 = x \\
	t(x,y,e_1',e_2') &= (x\wedge 0) \vee (y\wedge 1) = 0 \vee y = y
\end{align*}
\epf

While \textbf{Lat} has a diagonalising term, it is not coextensive. Note however, that the category \textbf{DLat} of distributive lattices is a subvariety of $\CSR$ closed under products, so \textbf{DLat} is coextensive. This example helps to illustrate how close any variety with a diagonalising term is to being co-extensive. In fact, having a diagonalising term is sufficient to imply the following properties:

\proposition\label{proposition_diag} Let $\C$ be a variety with a diagonalising term. Then the following conditions hold:
\begin{enumerate}
	\item $\C$ has co-disjoint products
	\item $\C$ is left co-extensive
	\item Any reflexive homomorphic relation $R$ on some product $A \times B \in \C$ is of the form $R_A \times R_B$ for some reflexive homomorphic relations $R_A$ on $A$ and $R_B$ on $B$.
	\item Any congruence $E$ on some product $A \times B \in \C$ is of the form $E_A \times E_B$ for some congruences $E_A$ on $A$ and $E_B$ on $B$.
	\item $\C$ satisfies the coextensivity condition for any surjective homomorphism $q$.
\end{enumerate}
\endthm

\pf Condition 1 holds since $k\geq1$, and so $e_k\in\Fo$ implies $\Fo$ is not empty. Condition 1 also follows from condition 2, which holds by Theorem 3 of \cite{DB} since we can fix any constant as $0$ and have the equalities:
	\begin{align*}
	    t(x, 0, e_1, \dots, e_k) = x \\
	    t(x, 0, e_1', \dots, e_k') = 0
	\end{align*}
Lastly, conditions 3,4, and 5 are simply proposition 2.3.
\epf

So, despite being a fairly simple condition, varieties with diagonalising terms are fairly close to being coextensive. Of course, unlike coextensivity, having a diagonalising term is a condition that only makes sense for varieties of algebras. Consider condition (3) of the above proposition. Any reflexive homomorphic relation $R$ on some product $A \times B \in \C$ is of the form $R_A \times R_B$ for some reflexive homomorphic relations $R_A$ on $A$ and $R_B$ on $B$. This condition can be made categorical in the following way:

\definition\label{definition_relation} Let $\C$ be a category. An internal relation from an object $A$ to an object $B$ is an object $R$ equipped with two morphisms $d_0: R \to A$ and $d_1: R \to B$ that are jointly monic. This is to say that given any object $C$ and morphisms $e_0: A \to C$ and $e_1: B \to C$, if $e_0d_0 = e_1d_1$ then $e_0 = e_1$.
\endthm

When $\C$ has binary products, then an internal relation $R$ from $A$ to $B$ is simply a subobject of $A \times B$, which is to say that $R$ is equipped with a monomorphism $d: R \to A\times B$. The morphisms $d_0$ and $d_1$ can then be recovered as $d_0 = \pi_A \circ d$ and $d_1 = \pi_B \circ d$. 

\vspace{1em}

An internal relation $d:R \to A\times A$ from $A$ to $A$ is called a relation on $A$, and $R$ is said to be reflexive if there exists a morphism $r: A \to R$ such that $d \circ r = \Delta_A : A \to A \times A$, where $\Delta_A$ is the diagonal morphism. Note that both $r: A \to R$ and $\Delta_A: A \to A^2$ are objects in the slice category $A/\C$, while $d: R \to A^2$ is a monomorphism from $r$ to $d$ in $A/\C$. A reflexive relation on $A$ is therefore a subobject of $\Delta_A$ in $A/\C$, and we denote the full subcategory of $A/\C$ containing the reflexive relations on $A$ as \textbf{RRel}$_{A}$. Similarly, let $i_B: B \to A$ be a monomorphism, then we write $\Delta_{A,B} = \Delta_A \circ i_B$ and say that a relation $(R,d)$ on $A$ is reflexive on $B$ if there exists some morphism $r_B: B \to R$ such that $d \circ r_B = \Delta_{A,B}$. Then, $r_B$ is a subobject of $\Delta_{A,B}$ in $B/\C$, and we denote the full subcategory of $B/\C$ containing the relations on $A$ that are reflexive on $B$ as \textbf{RRel}$_{A,B}$.

\proposition\label{product_reflect_mono} Let $\C$ be a category in which all product projections are epimorphisms. Then for any morphisms $f, g \in \C$, $f \times g$ is a monomorphism if and only if both $f$ and $g$ are monomorphisms.
\endthm
	
\pf It is clear that in any category if $f: X \to A$ and $g: Y \to B$ are monomorphisms, then so is $f \times g :X \times Y \to A \times B$. On the other hand, if $f \times g$ is a monomorphism, then for any $h,h': Z \to X$ such that $fh = fh'$ we can construct the diagram:
 \begin{equation}
\bfig
\node 000(0,600)[Z]
\node 001(600,600)[Z\times Y]
\node 002(1200,600)[Y]
\node 00(0,0)[X]
\node 01(600,0)[X\times Y]
\node 02(1200,0)[Y]
\node 10(0,-600)[A]
\node 11(600,-600)[A\times B]
\node 12(1200,-600)[B]
\arrow[001`000;\pi_Z]
\arrow[001`002;\pi_Y]
\arrow/@<-.2em>/[000`00;h]
\arrow/@<-.2em>/[001`01;h\times 1_Y]
\arrow|r|/@<.2em>/[000`00;h']
\arrow|r|/@<.2em>/[001`01;h'\times 1_Y]
\arrow[002`02;1_Y]
\arrow[01`00;\pi_X]
\arrow[01`02;\pi_Y]
\arrow[00`10;f]
\arrow[01`11;f\times g]
\arrow[02`12;g]
\arrow[11`10;\pi_A]
\arrow[11`12;\pi_B]
\efig
\end{equation}
In which $(f\times g)\circ(h\times 1_Y) = fh \times g = fh' \times g = (f\times g)\circ(h'\times 1_Y)$, so the two parralel morphisms $h\times 1_Y$ and $h'\times 1_Y$ are equal, since $f \times g$ is a monomorphism. Therefore $h\pi_z = \pi_X\circ(h\times1_Y) = \pi_X\circ(h'\times1_Y) = h'\pi_Z$ and since $\pi_Z$ is an epimorphism we have $h=h'$. Therefore $f$ is a monomorphism and likewise so is $g$ by repeating the same argument.
\epf

This result allows us to give a categorical version of proposition 3.3 condition 3. 

\proposition\label{categorical_3.3.3} Let $\C$ be a coextensive category, and let $i_X: X \to A$ and $i_Y: Y \to B$ be monomorphisms in $\C$. Then the canonical functor $\times\colon X/\C \times Y/\C \to (X \times Y)/\C$ can be restricted to the functor $\times\colon \textbf{RRel}_{A,X} \times \textbf{RRel}_{B,Y} \to \textbf{RRel}_{A \times B,X \times Y}$, and this functor is also an equivalence.
\endthm
\pf $\C$ is coextensive, so  $\times\colon X/\C \times Y/\C \to (X \times Y)/\C$ is an equivalence. It immediately follows that its restriction to $\textbf{RRel}_{A,X} \times \textbf{RRel}_{A,Y}$ is fully faithful, and it remains to show that the restriction is essentially surjective. $\times$ itself is essentially surjective, so any $r: X\times Y \to R$ is isomorphic in $(X\times Y)/\C$ to some product $r_X \times r_Y: X\times Y \to R_X \times R_Y$. Note that $r \in \textbf{RRel}_{A\times B, X\times Y}$ if and only if there exists some monomorphism $d : R_X \times R_Y \to (A\times B)^2$ such that $dr = \Delta_{A\times B,X \times Y} = \Delta_{A\times B} \circ (i_X\times i_Y) \simeq \Delta_Ai_X \times \Delta_Bi_Y = \Delta_{A,X} \times \Delta_{B,Y}$. We have that $\times$ is fully faithful, so by the previous proposition this is true if and only if $d = d_X \times d_Y$ for some monomorphisms $d_A: R_X \to A^2$ and $d_B : R_Y \to B^2$ such that $d_Ar_X = \Delta_{A,X}$ and $d_Br_Y = \Delta_{B,Y}$, which is to say that $(r_X,r_Y) \in \textbf{RRel}_{A,X} \times \textbf{RRel}_{B,Y}$. Therefore the restriction is both fully faithful and essentially surjective, so it is an equivalence as desired.
\epf

In the case where $X=A$ and $Y=B$ we have:

\corollary\label{cor_reflexive} Let $\C$ be a coextensive category. Then for any $X,Y\in \C$ the canonical functor $\times\colon \textbf{RRel}_{X} \times \textbf{RRel}_{Y} \to \textbf{RRel}_{X \times Y}$ is an equivalence.
\endthm

Note that a reflexive internal relation in a variety of universal algebras is exactly a reflexive homomorphic relation. Therefore the above result is simply a more general form of condition (3) of proposition 3.3. On the other hand, if a category $\C$ has an initial object $0$, then any object $X$ in $\C$ has a unique arrow $!_X : 0 \to X$, and in particular any relation on $A$ is an object of $\textbf{RRel}_{X,0}$, which we denote as $\textbf{Rel}_{A}$ and refer to as the category of relations on $A$. Note, however, that while the product of two reflexive homomorphic relations in a variety is simply a reflexive homomorphic relation, the product of two relations $R_X \times R_Y$ is always reflexive on pairs of constants. This is because for any two constants $e,e'$, both $R_X$ and $R_Y$ must contain the pairs $(e,e)$ and $(e,e')$, and therefore $((e,e'),(e,e')) \in R_X \times R_Y$. More generally, for any category $\C$, since $r_X \circ !_{R_X} = !_{X^2}$ and $r_Y \circ !_{R_Y} = !_{Y^2}$, we have $(r_X \times r_Y) \circ !_{R_X \times R_Y} = !_{X^2} \times !_{Y^2}$, so the product of two relations must be reflexive on some subproduct. Specifically, we say that a relation $R$ on some product $X \times Y$ is \textbf{subproduct reflexive} if and only if it is a subobject of $!_X \times !_Y$ in $(0\times 0)/\C$, and we denote by \textbf{SPRRel}$_{X,Y}$ the full subcategory $\textbf{RRel}_{X \times Y,0\times0}$ of $(0\times 0)/\C$ containing all such relations. It then follows as a corollary of the previous proposition:

\corollary\label{cor_subproduct_reflexive} Let $\C$ be a coextensive category with an initial object. Then for any $X,Y\in \C$ the canonical functor $\times\colon \textbf{Rel}_{X} \times \textbf{Rel}_{Y} \to \textbf{SPRRel}_{X,Y}$ is an equivalence.
\endthm

And if $\C$ is a variety, then this corollary is equivalent to the existence of a diagonalising term as follows.

\proposition\label{categorical_diagonalising} Let $\C$ be a variety of universal algebras, then $\C$ has a diagonalising term if and only if the canonical functor $\times\colon 0/\C \times 0/\C \to (0 \times 0)/\C$ can be restricted to the functor $\times'\colon \textbf{Rel}_{X} \times \textbf{Rel}_{Y} \to \textbf{SPRRel}_{X, Y}$, and this functor is also an equivalence.
\endthm
\pf Note that $\times'$ is essentially surjective if and only if, for any subproduct reflexive relation $R$ on $X\times Y$, there exist some relations $R_x$ on $X$ and $R_Y$ on $Y$ and morphisms $p_1: R \to R_X$, $p_2: R \to R_Y$ such that the middle row of the following diagram is a product diagram. 

\begin{equation}
\bfig
\node 03(1800,0)[0]
\node 04(2400,0)[0\times 0]
\node 05(3000,0)[0]
\node 13(1800,-600)[R_X]
\node 14(2400,-600)[R]
\node 15(3000,-600)[R_Y]
\node 23(1800,-1200)[X^2]
\node 24(2400,-1200)[(X\times Y)^2]
\node 25(3000,-1200)[Y^2]
\arrow[04`03;\pi_1]
\arrow[04`05;\pi_2]
\arrow[03`13;!_{R_X}]
\arrow[04`14;h]
\arrow[05`15;!_{R_Y}]
\arrow[14`13;p_1]
\arrow[14`15;p_2]
\arrow[13`23;r_X]
\arrow[14`24;r_{X\times Y}]
\arrow[15`25;r_Y]
\arrow[24`23;\pi_X]
\arrow[24`25;\pi_Y]
\efig
\end{equation}

Since $\C$ is a variety, we can calculate $R_X$ and $R_Y$ as before as:
\begin{align*}
	R_X &= \{(a,c) \in X^2 | ((a,b),(c,d)) \in R \text{ for some } b,d \in Y\} \\
	R_Y &= \{(b,d) \in X^2 | ((a,b),(c,d)) \in R \text{ for some } a,c \in X\}
\end{align*}
Therefore, their product $R_X \times R_Y$ is of the form:
\begin{align*}
	((a,b),(c,d)) \in R_X \times R_Y &\iff (a,c) \in R_X \text{ and } (b,d) \in R_Y \\
	&\iff ((a,b'),(c,d')) \in R \text{ and } ((a',b),(c',d)) \in R
\end{align*}
Let $X = F(\{a,a',c,c'\})$, $Y = F(\{b,b',d,d'\})$ and $R$ be the relation generated by $(a,b')\simeq_R(c,d'), (a',b)\simeq_R(c',d)$, and $(e,e') \simeq (e,e')$ for all constants $e,e' \in \Fo$. Then $R = R_X \times R_Y$ if and only if $(a,b) \simeq_R (c,d)$ if and only if: $$((a,b),(c,d)) = t(((a,b'),(c,d')),((a',b),(c',d)),((e_1,e_1'),(e_1,e_1')),\dots,((e_k,e_k'),(e_k,e_k')) \in R$$
for some term $k \in \N$ and $k+2$-ary term $t$. However this is simply to say that $\C$ satisfies the identities:
\begin{align*}
    t(x, y, e_1, \dots, e_k) &= x \\
    t(x, y, e_1', \dots, e_k') &= y
\end{align*}
or in other words, $t$ is a diagonalising term. On the other hand, if $\C$ has a diagonalising term then for any $R$, $R \simeq R_X \times R_Y$ where $R_X$ and $R_Y$ are defined as above. Therefore $\times'$ is essentially surjective. Furthermore, since $\C$ has a diagonalising term, $\C$ is left coextensive. This is to say that $\times$ is fully faithful, and therefore so is its restriction $\times'$. Therefore $\times'$ is an equivalence if and only if $\C$ has a diagonalising term as desired.
\epf

\refs

\bibitem [1]{CLW} [1] A. Carboni, S. Lack, and R. F. C. Walters, Introduction to extensive and distributive categories, Journal of Pure and Applied Algebra, 84(2), 1993, 145-158.

\bibitem [2]{DB} [2] D. Broodryk, Characterization of left coextensive varieties of universal algebras, Theory and Applications of Categories, Vol. 34, No. 32, 2019, 1036-1038.

\bibitem [3]{COEX} [3] D. Broodryk, Characterization of coextensive varieties of universal algebras, arXiv:2008.03474 [math.CT].

\endrefs

\end{document}